\input AHTOH-E.STY
\hfuzz7pt
\def\Epi{{\rm Epi}}
\def\Mono{{\rm Mono}}

\UDC{
512.542        
+ 512.543.72   
}


\title{%
On the number of epi-, mono-, and homomorphisms of groups
}
\author{%
Elena K. Brusyanskaya
\qquad
\qquad
Anton A. Klyachko
}
\address{
Faculty of Mechanics and Mathematics of Moscow State University
\\
Moscow 119991, Leninskie gory, MSU.
\\
Moscow Center for Fundamental and Applied Mathematics 
\\
\strut
ebrusianskaia@gmail.com
\quad
klyachko@mech.math.msu.su
}

\grants{\RFBR19-01-00591}%

\abstract{%
It is known that
the number of homomorphisms from a group $F$ to a group~$G$
is divisible by the greatest common divisor of the
order of~$G$ and the exponent of~$F/[F,F]$.
We investigate the
number of homomorphisms
satisfying some natural conditions such as injectivity or
surjectivity.
The simplest nontrivial corollary of our results
is the following fact:
{\it in any finite group, the
number of generating pairs $(x,y)$ such that
$x^3=1=y^5$,
is a multiple of the greatest common divisor
of 15 and the order of the group
$[G,G]\cdot\{g^{15}\;|\;g\in G\}$.
}
}

\s 0.
Introduction

We were inspired by the three classical results
on the divisibility in groups:
the theorems of Frobenius~(1895), Solomon~(1969), and Iwasaki~(1985).

\proclaim Frobenius theorem {\rm[Frob95] (see also [And16])}.
The number of solutions to an
equation $x^n=1$ in a finite group~$G$
is divisible by $\GCD(|G|,n)$
for any integer $n$.

\proclaim Solomon theorem \rm [Solo69].
In any
group,
the number of solutions to a finite system of coefficient-free equations
is divisible by the order of the group if
the number of equations is fewer than
the number of
unknowns.
\newline
{\rm In other words,} the number of homomorphisms
$\pres<x_1,\dots,x_m|w_1=\dots=w_n=1>\to G$
is divisible by $|G|$ if $m>n$.

\proclaim Iwasaki theorem \rm [Iwa82].
For any integer $n$, the number of elements
of a finite group $G$ whose $n$-th powers
lie in a subgroup $A\subseteq G$ is divisible by $|A|$.

\noindent
These theorems were generalised in several directions,
see, e.g.,
[Frob03],
[Hall36],
[Kula38],
[Sehg62],
[Isaa70],
[BrTh88],
[Yosh93],
[Stru95],
[AsTa01],
[SaAs07],
[AmV11],
[GRV12],
[ACNT13]
[KM14],
[KM17],
[BKV19],
[KR20],
and literature cited therein.
For example, in [GRV12],
the following generalisation of the Solomon theorem was obtained.

\proclaim Gordon--Rodriguez-Villegas theorem {\rm [GRV12]}.
The number of homomorphisms $F\to G$ is divisible by the order of
the group $G$
if $F$ is a finitely generated group whose
commutator subgroup is of infinite index.

\noindent
Later, it turned out that the three classical results are related to each
other:
\-
in [KM17], it was obtained a general fact, called
\emph{Theorem KM} here, that includes
as special cases
Solomon's and Iwasaki's theorem (as well as their generalisations
and some other interesting facts);
\-
in [BKV19], it was shown that all three classical theorem (as well as
their
generalisations including Theorem KM) are special cases
of a very general
fact, which we call \emph{Theorem BKV} (see the following
section).

\enditem
The authors of [KM17] applied Theorem KM
to obtain
the
following result on the divisibility of the number of
homomorphisms satisfying injectivity or
surjectivity conditions.
Suppose that $F\supseteq W$ and $G\supseteq A$ are groups and
$$
\eqalign{
\Hom(F,W;G,A)&=\{\phi\:F\to G\;|\;\phi(W)\subseteq A\},
\quad
\Epi(F,W;G,A)=\{\phi\:F\to G\;|\;\phi(W)=A\},
\cr
\Mono(F,W;G,A)&=
\{\phi\:F\to G\;|\;\phi(W)\subseteq A\
\hbox{and the restriction of $\phi$ to $W$ is injective}
\}.
}
$$
\Theorem on epi-, mono-, and homomorphisms {\rm[KM17]}.
Let $A$ be a subgroup of a group $G$, and let $W$ be a subgroup
of a finitely
generated group $F$ whose
commutator subgroup $F'$ is of
infinite index.
Then
\item{\rm a)}
$|\Hom(F,W;G,A)|$, $|\Epi(F,W;G,A)|$, and $|\Mono(F,W;G,A)|$
are divisible by the order of the normaliser $N(A)$ of $A$
if the index $|F:F'W|$ is infinite\;
\item{\rm b)}
$|\Hom(F,W;G,A)|$ is divisible by $|A|$\;
\item{\rm c)}
$|\Epi(F,W;G,A)|$ is divisible by $|A'|$.

\noindent
The purpose of this paper is
to
add ``Frobeniusness" to this theorem, i.e. to get rid of the
conditions $|F:F'|=\infty$ and $|F:F'W|=\infty$.
The answer turns out to be 
expected For Assertions a) and b), 
less obvious for c); in addition, a new fact d) arises.

\proclaim ``Frobeniusian" theorem on epi-, mono-, and homomorphisms.
Let $A$ be a subgroup of a group $G$, and let $W$ be a subgroup
of a finitely
generated group $F$.
Then
\item{\rm a)}
$|\Hom(F,W;G,A)|$, $|\Epi(F,W;G,A)|$, and $|\Mono(F,W;G,A)|$
are
divisible by
$\GCD\Bigl(N(A),\exp\bigl(F/(F'W)\bigr)\Bigr)$\;
\item{\rm b)}
$|\Hom(F,W;G,A)|$
is divisible by
$\GCD\bigl(A,\exp(F/F')\bigr)$\;
\item{\rm c)}
$|\Epi(F,W;G,A)|$
is divisible by
$\GCD\bigl(A'A^{\exp(F/F')},\exp(F/F')\bigr)$\;
\item{\rm d)}
$|\Mono(F,W;G,A)|$
is divisible by
$\GCD\biggl(A,\exp\Bigl(F/\bigl(F'Z(W)\bigr)\Bigr)\biggr)$.

\noindent
This is an over-simplified
version of Theorem 1 or, to be more precise, of its
corollary, see Section~3.
Both assertions about $|\Hom(F,W;G,A)|$ are not new --- these facts
established in~[BKV19] (in somewhat another language) 
are included here for the sake of completeness.

Note that we do not assume that the group $G$ is
finite. We follow notation of [BKV19]: the
\emph{greatest common divisor $\GCD(G,n)$} of a group $G$ and an integer
$n$ is the least common multiple of the orders of subgroups of $G$
dividing $n$; the divisibility is always understood in the sense
of cardinal arithmetics: each infinite cardinal is divisible by all
smaller nonzero cardinals (and zero is surely divisible by all
cardinals and divides only zero). This means that
$\GCD(G,0)=|G|$ for any group~$G$; and, e.g., 
$\GCD(\SL_2(\Z),2020)=2$.
Though, the reader will not lose very much by
assuming all groups to be finite; in this case,
$\GCD(G,n)=\GCD(|G|,n)$ by Sylow's theorem (and because a finite
$p$-group contains subgroups of all possible orders).

\noindent
Assertion b) of this theorem contains all
three classical results:
\-
Frobenius's theorem
(take a cyclic group as $F=W$
and put $A=G$),
\-
Solomon's and even the Gordon--Rodriguez-Villegas theorem
(take a finitely generated group whose commutator subgroup is
of infinite index as $F=W$ and put $A=G$),
\-
and Iwasaki's theorem (put $F=\Z\supseteq n\Z=W$).

\enditem
If, in Assertion c), we take a
free product cyclic groups
as $F=W$
and put $A=G$, then
we obtain the following fact.

\proclaim{Corollary on generating tuples}.
For any group $G$, and any $k_i\in\Z$ the number
of tuples
$(g_1,\dots,g_n)$ of elements of $G$ such that
$\gp{g_1,\dots,g_n}=G$ and $g_i^{k_i}=1$
is divisible by
$\GCD\(G'\cdot G^{\LCM(k_1,\dots,\,k_n)},\LCM(k_1,\dots,\,k_n)\)$.
\rm(Henceforth, $G^m\:=\gp{\{g^m\;|\;g\in G\}}$.)

The reader can guess that the key to our generalisation
of the theorem on
epi-, mono-, and homomorphisms is the use instead of Theorem KM
its ``Frobeniusian analogue", i.e. Theorem BKV.
This is true, but
actually, we generalise Theorem BKV itself,
see the main theorem in the following section and its proof
in Section~2.

{\noindent \bf Notation and conventions}
we use are mainly standard. Note only that, if
$k\in \Z$ and $x,y$ are elements of a group, then $x^y$, $x^{ky}$ and
$x^{-y}$ denote $y^{-1}xy$, $y^{-1}x^ky$ and $y^{-1}x^{-1}y$,
respectively.
The commutator subgroup of a group~$G$ is denoted by $G'$
or $[G,G]$; the centre of a group $G$ is denoted as~$Z(G)$.
The subgroup of $G$ generated by the $n$th powers of
all elements is denoted by~$G^n$.
The symbol $|X|$ denotes the cardinality of a set $X$.
If $X$ is a
subset of a group, then $\gp X$,
$C(X)$ and~$N(X)$
mean
the subgroup generated by~$X$,
centraliser of~$X$,
and the normaliser of~$X$,
respectively. The index of a subgroup $H$ of a
group~$G$ is denoted by $|G:H|$.  The letter~$\Z$ denotes
the set of integers.
$\GCD$ and $\LCM$~ are the greatest common divisor and the least
common multiple. The symbol $\exp(G)$ stands for the period
(exponent) of a
group~$G$ if the period is finite; we assume $\exp(G)=0$ if the period
is infinite. Let us repeat that the finiteness of groups are not
assumed by default; the divisibility is always understood in the sense
of cardinal arithmetics (an infinite cardinal is divisible by all
nonzero cardinals not exceeding it); and~%
$
\GCD(G,n)\:=\LCM\(\bigl\{|H|\;\bigm|\;
\hbox{$H$ is a subgroup of $G$ and $|H|$ divides $n$}\bigr\}\).
$

\s 1.
Main theorem

A group $F$ equipped with an epimorphism $F\to\Z_n\:=\Z/n\Z$
(where $n\in\Z$)
is called an
\emph{$n$-indexed} group [BKV19].
This epimorphism $F\to\Z_n$ is called
\emph{degree} and denote $\deg$.
Thus, for any element~$f$ of an
indexed group $F$, an element~$\deg f\in\Z_n$ is assigned;
$F$ contains elements of all degrees, and
${\deg(fg)=\deg f+\deg g}$ for any $f,g\in F$.

Suppose that $\phi\:F\to G$ is a homomorphism from an $n$-indexed
group $F$ to a group $G$ containing a subgroup
$H$. The subgroup
$
H_\phi=\bigcap\limits_{f\in F}\!H^{\phi(f)}\cap C\bigl(\phi(\ker\deg)\bigr)
$
is called the
\emph{$\phi$-core} of~$H$ [KM17]. In other words,
the $\phi$-core $H_\phi$ of $H$ consists of
elements~$h\in H$ such that $h^{\phi(f)}\in H$ for all $f$, and
$h^{\phi(f)}=h$ if $\deg f=0$.

\Theorem BKV {\rm[BKV19]}.
Suppose that an
integer $n$
is divisible by the order of a subgroup
$H$ of a group $G$,
and a
set $\Phi$ of homomorphisms from
an $n$-indexed group $F$ to $G$
satisfies the following
conditions.
\item{\rm I.}
$\Phi$ is invariant with respect to conjugation
by elements of $H$\:
if $h\in H$ and $\phi\in\Phi$, then the homomorphism
$\psi\:f\mapsto\phi(f)^h$ belongs to $\Phi$.
\item{\rm II.}
For any $\phi\in\Phi$ and any element $h$ of the $\phi$-core
$H_\phi$ of $H$, the homomorphism $\psi$ defined by
$$
\psi(f)=
\cases{
\phi(f)& for all elements $f\in F$ of degree zero;
\cr
\phi(f)h& for
some element
$f\in F$ of degree one
\small(and, hence, for all elements of degree one),
\cr
}
$$
lies in $\Phi$.
\enditem
Then $|\Phi|$ is divisible by
$|H|$.

\goodbreak

\noindent
Note, that
\-
the mapping $\psi$ from Condition I is a homomorphism for
any $h\in G$;
and the formula for~$\psi$ from Condition II defines a homomorphism
for any $h\in H_\phi$ (as explained in [BKV19]);
thus Conditions~I~and~II
only require these homomorphisms to lie in $\Phi$;
\-
according to (a simple) Lemma~3 from~[BKV19],
{\sl
in Condition {\rm II},
$\psi(f)\in\phi(f)H_\phi$ for all
$f\in F$;
}
\-
the condition ``$n$
is divisible by the order of $H$"
can be omitted, but then the conclusion of the theorem should be:
``$|\Phi|$ is divisible by
$\GCD(H,n)$" (instead of ``$|\Phi|$ is divisible by $|H|$");
this follows immediately from the definition of the greatest common
divisor  of a
group and an integer (see Introduction), because,
if Conditions I and II hold for $H$, then they hold also
for any subgroup of~$H$;
\-
Theorem KM (mentioned in
Introduction) is exactly Theorem BKV with $n=0$.

\proclaim{Main theorem}.
Suppose that
$F$ is an $n$-indexed group,
$H$ is a subgroup group $G$,
$k$ is a positive integer,
and
$\Phi$ is a set of homomorphisms from $F$ to $G$
satisfying the following conditions\:
\item{\rm(i)}
for all
$\phi\in\Phi$
and
$h\in H$,
the homomorphism
$\psi\:f\mapsto\phi(f)^h$ lies in $\Phi$\;
\item{\rm(ii)}
for each $\phi\in\Phi$, the
$\phi$-core $H_\phi$
of $H$
contains
a subgroup $H_{\phi,k}$
such that
\itemitem{-}
$H_\phi\supseteq H_{\phi,k}\nin\gp{H_\phi\cup\phi(F)}$;
\itemitem{-}
$|H_\phi/H_{\phi,k}|$ divides $k$;
\itemitem{-}
if $\phi\in\Phi$ and $\psi\:F\to G$ is a
homomorphism coinciding with
$\phi$ on elements of degree zero and
such that
$\psi(w)\in\phi(w)H_{\phi,k}$ for all
elements~$w\in F$ whose degrees are multiples of $k$
\(i.e. $\deg(w)\in k\Z_n$\),
then $\psi\in\Phi$.
\enditem
Then $|\Phi|$ is divisible by $\GCD(H,n)$.

For $k=1$ and $H_{\phi,k}=H_\phi$,
this fact is effectively Theorem BKV.

\s 2.
Proof of the main theorem

We can assume that $|H|$ divides $n$
(by the definition of the greatest common divisor of a group and an
integer, because Conditions (i) and (ii) remains satisfied when we
replace $H$ by
its subgroup). Now, it suffices to show that Conditions~I~and~II
of Theorem BKV hold for these $F$, $G$, $H$, and $\Phi$. Condition~I
holds obviously by virtue of Condition~(i).

Let us verify Condition II.
Suppose that $\phi\in\Phi$, an element $f_1\in F$ has degree
one, $\phi(f_1)=g$, and $h\in H_\phi$. We have to show that the
homomorphism
$\psi\:F\to G$ coinciding with $\phi$ on elements of degree zero and
mapping $f_1$ to $gh$ belongs to $\Phi$.
Each element $w\in F$ whose degree is divisible by $k$
can be written in the form
$w=f_0f_1^{ki}$ for some $i\in\Z$ and $f_0\in\ker\deg$.
Then
$$
\psi(w)=\psi(f_0f_1^{ki})=\psi(f_0)(gh)^{ki}=\phi(f_0)(gh)^{ki}
\qbox{(because $\phi$ and $\psi$ coincide on $\ker\deg$)}.
$$
The subgroup $H_\phi$ is normal in $\gp{H_\phi,g}$
by the definition of $\phi$-core $H_\phi$.

\proclaim Brauer lemma {\rm[Bra69] (see also [BKV19])}.
If $U$ is a finite
normal subgroup of a group $V$,
then, for all $v\in V$ and $u\in U$,
the elements $v^{|U|}$ and $(vu)^{|U|}$ are conjugate
by an element of $U$.

Applying the Brauer Lemma
to the normal
subgroup $H_\phi/H_{\phi,k}$ of
$\gp{g, H_\phi}/H_{\phi,k}$,
we obtain that
$(gh)^{ki}\in g^{kih'}H_{\phi,k}$ for some $h'\in H_\phi$,
which depends on neither $i$ nor $w$, and is determined
by the homomorphisms $\phi$ and $\psi$ only.
Therefore,
$$
\psi(w)
=
\phi(f_0)(gh)^{ki}
\in
\phi(f_0)g^{kih'}H_{\phi,k}
\=^1
\(\phi(f_0)g^{ki}\)^{h'}H_{\phi,k}
\=^2
\(\phi(f_0f_1^{ki})\)^{h'}H_{\phi,k}
\=^3
\bigl(\phi(w)\bigr)^{h'}H_{\phi,k},
\hbox{ where}
$$
\-
$\=^1$
holds because
$h'\in H_\phi$
commutes with the image $\phi(f_0)$ of the element $f_0$
of degree zero
by the definition of the $\phi$-core $H_\phi$;
\-
$\=^2$
follows from the definition of $g$;
\-
$\=^3$
follows from the definition of $w$.

\enditem
The homomorphism $f\mapsto \bigl(\phi(f)\bigr)^{h'}$
lies in $\Phi$ by Condition (i), and, therefore,
$\psi\in\Phi$ by Condition (ii).
This completes the proof.

\s 3.
What is the number of epi-, mono-, and homomorphisms divisible by?

Suppose that $\Phi$ is a set of homomorphisms
from an $n$-indexed group $F$ to a group $G$ containing
subgroups $B$ and $H$.
We say that the subgroup $H\subseteq G$ is \emph{$(B,k,\Phi)$-smooth} if,
for all $\phi\in\Phi$,
the subgroup
$H_\phi\cap B$
contains a normal in $\gp{H_\phi,\phi(F)}$
subgroup $\^B$
\(depending on $\phi$\)
such that
$|H_\phi/\^B|$
divides $k$.

The following lemma consists of some obvious examples of smooth subgroups.

\Lemma on smooth subgroups.
The following subgroups of a group $G$ are $(B,k,\Phi)$-smooth\:
\item{\rm1)}
any subgroup contained in $B$\;
\item{\rm2)}
any subgroup of order dividing $k$\;
\item{\rm3)}
any subgroup $H$ such that $|H:H\cap B|$
divides $k$, if $B\nin G$.

\Proof
We can take the following subgroups as $\^B$:
1) $H_\phi$;
\quad 
2) $\1$;
\quad 
3) $H\cap B$.

\Th 1.
Suppose that a group~$G$ contains a subgroup $A$,
an $n$-indexed group~$F$ contains a subgroup $W$
such that
$\deg(W)=k\Z_n$, where $k$ divides $n$,
and
$$
\eqalign{
\Hom(F,W;G,A)&=\{\phi\:F\to G\;|\;\phi(W)\subseteq A\},
\quad
\Epi(F,W;G,A)=\{\phi\:F\to G\;|\;\phi(W)=A\},
\cr
\Mono(F,W;G,A)&=
\{\phi\:F\to G\;|\;\phi(W)\subseteq A\
\hbox{and the restriction of $\phi$ to $W$ is injective}
\}.
}
$$
Then $\GCD(H,n)$ divides
\item{\rm a)}
$|\Hom(F,W;G,A)|$
for any $\bigl(A,k,\Hom(F,W;G,A)\bigr)$-smooth subgroup
$H\subseteq N(A)$ of
$G$\;
\item{\rm b)}
$|\Epi(F,W;G,A)|$
for any
$\bigl(A'A^n,k,\Epi(F,W;G,A)\bigr)$-smooth subgroup
$H\subseteq N(A)$ of
$G$, where $A^n\:=\gp{\{a^n\;|\;a\in A\}}$\;
\item{\rm c)}
$|\Mono(F,W;G,A)|$
for any $\bigl(A,k,\Mono(F,W;G,A)\bigr)$-smooth
subgroup $H\subseteq N(A)$ of $G$ if the indexation
of $F$ is such that $\deg(w)=0$
for each central \(in $W$\) element $w\in W$
such
that $w^n=1$.

\goodbreak
\noindent
(Note parenthetically that
the subgroup $A'A^n$ from Assertion~b)
and the subgroup $\{w\in Z(W)\;|\;w^n=1\}$ from Assertion c)
are the verbal subgroup of $A$ and marginal
subgroup of $W$ corresponding to the variety of abelian groups
of exponent~$n$.)

\Proof
It suffices to
verify that Conditions (i) and (ii)
of the main theorem hold for given
$F$, $G$, $H$, $k$, $\Phi$ and $H_{\phi,k}=\^B$
(here, $\^B$ is from definition smooth subgroup $B$, where
$B=A$
in assertions a) and c), and $B=A'A^n$ in assertion b)).
The first two items of Condition (ii) are satisfied
by the definition of smoothness; we need to verify only the last
part of
Condition~(ii).

\item{a)} $B=A$ and $\Phi=\Hom(F,W;G,A)$.
Condition (i)
holds  obviously, because $H\subseteq N(A)$. Condition~(ii) also holds,
because. for all $w\in W$, we have
$\psi(w)\in\phi(w)\^B\subseteq\phi(w)A=A$, i.e. $\psi\in\Phi$ as
required.

\item{b)}
$B=A'A^n$ and $\Phi=\Epi(F,W;G,A)$.
Condition (i) holds obviously by the same reason: $H\subseteq N(A)$.
Condition~(ii) also holds:
$
A\=^1\phi(W)\subseteqq^2\psi(W)A'A^n\=^3\psi(W)\phi(W'W^n)\=^4
\psi(W)\psi(W'W^n)=
\psi(W),
$
where
\itemitem{$\=^1$}
holds by the definition of $\Phi\ni\phi$;
\itemitem{$\subseteqq^2$}
holds by the definition of $\psi$ from Condition (ii) for
$B=A'A^n\supseteq\^B=H_{\phi,k}$;
\itemitem{$\=^3$}
follows from $\=^1$;
\itemitem{$\=^4$}
holds because $\deg(W'W^n)=\0$, and
homomorphisms
$\psi$ and $\phi$ from Condition (ii)
coincide on elements of degree zero.

So, we obtain that $A=\psi(W)$, i.e. $\psi\in\Phi$ as required.

\item{c)}
$B=A$ and $\Phi=\Mono(F,W;G,A)$.
Condition (i) holds obviously by the same reason: $H\subseteq N(A)$.
Let us show that (ii) also holds.
First,
$\psi(W)\subseteq\phi(W)A=A$.
It remains to show that $\ker\psi\cap W=\1$.
Take $w\in\ker\psi\cap W$. Then
\itemitem{-}
for each $w'\in W$, we have
$1=\psi([w,w'])=\phi([w,w'])$ (because commutators have degree zero,
and $\phi$ and $\psi$ coincide on elements of degree zero), therefore,
$[w,w']=1$ (because $\phi$ is injective on $W$),
i.e. $w\in Z(W)$;
\itemitem{-}
similarly, we obtain
$1=\psi(w^n)=\phi(w^n)$ (because $\deg(w^n)=n\deg(w)=0$,
and  $\phi$ and $\psi$
coincide on elements of degree zero), therefore,
$w^n=1$ (because $\phi$ is
injective on $W$).

\item{}
We obtain that $w$ is a central element of the group $W$ and $w^n=1$;
such elements have degree zero by the condition. Therefore,
$\phi(w)=\psi(w)=1$, i.e. $w=1$, since
$\phi\in\Mono(F,W:G,A)$. Therefore, $\ker\psi\cap W=\1$ as
required.

\Corollary.
In conditions of Theorem 1,
$|\Hom(F,W;G,A)|$,
$|\Epi(F,W;G,A)|$,
and $|\Mono(F,W;G,A)|$
are divisible by $\GCD\bigl(k,\,N(A)\bigr)$.
Moreover,
\item{\rm a)}
$|\Hom(F,W;G,A)|$
is divisible
\itemitem{-}
by
$\GCD\bigl(n,\;A\bigr)$,
\itemitem{-}
and also by
$
\GCD\Bigl(n,\;
|A|\cdot\GCD\bigl(k,\;|G/A|\bigr)\Bigr)
=
\GCD(n,\ |G|,\ k\cdot|A|)
$
if $A\nin G$ and $G$ finite\;
\item{\rm b)}
$|\Epi(F,W;G,A)|$
is divisible
\itemitem{-}
by
$\GCD\bigl(n,\;A'A^n\bigr)$
and even by $\GCD\bigl(n,\;H\bigr)$,
where $H$ is any subgroup of $N(A)$ such
that
\newline
$|(C(A'A^n)\cap H:Z(A'A^n)\cap H|$ divides $k$,
\itemitem{-}
and also by
$
\GCD\Bigl(n,\;
|A'A^n|\cdot\GCD\bigl(k,\;|G/(A'A^n)|\bigr)\Bigr)
=
\GCD(n,\ |G|,\ k\cdot|A'A^n|)
$,
if $A\nin G$ and $G$ is finite\;
\item{\rm c)}
if $\deg\bigl(\{w\in Z(W)\;|\;w^n=1\}\bigr)=\0$, then
$|\Mono(F,W;G,A)|$
is divisible
\itemitem{-}
by
$\GCD\bigl(n,\;A\bigr)$,
\itemitem{-}
and also by
$
\GCD\Bigl(n,\,
|A|\cdot\GCD\bigl(k,\,|G/A|\bigr)\Bigr)
=
\GCD(n,\,|G|,\,k\cdot|A|)
$
if $A\nin G$ and $G$ is finite.
\enditem

\Proof
The first assertion (on the divisibility by $\GCD\bigl(k,\,N(A)\bigr)$)
follows immediately from Theorem 1 and Assertion~2) of
the lemma on smooth
subgroups.

The remaining assertion of this corollary follow also from Theorem 1
and a suitable fact about smooth subgroups.
\item{a)}
The divisibility by $\GCD\bigl(n,\;A\bigr)$
follows immediately from Assertion 1) of
the lemma on smooth subgroups.
The divisibility by $\GCD(n,\;|G|,\;k\cdot|A|)$ follows from
Assertion 3) of the lemma on smooth subgroups.
Indeed, In Theorem~1, take
$H$ to be a $p$-subgroup of $G$
whose
order is the maximal power $p^i$ of a prime~$p$
dividing $\GCD(n,k|A|,|G|)$, and choose $H$
\itemitem{-}
inside $A$ if $p^i$
divides $|A|$;
\itemitem{-}
containing the Sylow $p$-subgroup of $A$ otherwise.

\goodbreak

\item{}
This subgroup is
\newline
$\bigl(A,k,\Hom(F,W;G,A)\bigr)$-smooth
by the lemma on smooth subgroups.
Therefore, $|H|$ divides $\Hom(F,W;G,A)$ by
Theorem~1. Doing so for each prime $p$, we
obtain the required divisibility.

\goodbreak

\item{b)}
The second assertion of b) is proved
by the very same way as the second assertion of a).
To prove the first assertion of b)
using
Theorem~1, it suffices to show that $H$ is
$\bigl(A'A^n,k,\Epi(F,W;G,A)\bigr)$-smooth,
i.e.,
for all $\phi\in\Epi(F,W;G,A)$,
the group
$H_\phi\cap(A'A^n)$
contains a normal in $\gp{H_\phi,\phi(F)}$
subgroup $\^B$
such that
$|H_\phi/\^B|$
divides $k$.
We can take $\^B=Z(A'A^n)\cap H_\phi$.
Indeed,
$$
H_\phi
\subseteqq^1
C\bigl(\phi(\ker\deg)\bigr)
\subseteqq^2
C\bigl(\phi(W'W^n)\bigr)
\=^3
C\bigl(A'A^n\bigr),
\qbox{where}
$$
{$\subseteqq^1$}
follows from the definition of the $\phi$-core $H_\phi$,
{$\subseteqq^2$}
holds as $\deg(W'W^n)=\0$,
and
{$\=^3$}
holds as $\phi(W)=A$.

\item{}
Therefore, by the Lagrange theorem,
$\left|H_\phi/\bigl(Z(A'A^n)\cap H_\phi\bigr)\right|$
divides
$\left|\bigl(C(A'A^n)\cap H\bigr)/
\bigl(Z(A'A^n)\cap H\bigr)\right|$
that divides $k$ by the condition.

\item{c)}
Here, the argument is completely similar to a).

\baselineskip 10.3pt

\References

[AmV11]
A. Amit, U. Vishne,
Characters and solutions to equations in finite groups,
J. Algebra Appl., 10:4 (2011), 675-686.

[And16]
R. Andreev,
A translation of
``Verallgemeinerung des Sylow'schen Satzes"
by F. G. Frobenius.
\newline
arXiv:1608.08813.

[ACNT13]
T. Asai, N. Chigira, T. Niwasaki, Yu. Takegahara,
On a theorem of P. Hall,
Journal of Group Theory, 16:1 (2013), 69-80.

[AsTa01]
T. Asai, Yu. Takegahara,
$|\Hom(A,G)|$, IV,
J. Algebra, 246 (2001), 543-563.


[Bra69]
R. Brauer,
On A Theorem of Frobenius,
The American Mathematical Monthly, 76:1 (1969), 12-15.

[BrTh88]
K. Brown, J. Th\'evenaz,
A generalization of Sylow's third theorem,
J. Algebra, 115 (1988), 414-430.

[BKV19]
E. K. Brusyanskaya, A. A. Klyachko, A.V. Vasil'ev,
What do Frobenius's, Solomon's, and Iwasaki's theorems on divisibility
in groups have in common?,
Pacific Journal of Mathematics, 302:2 (2019), 437-452.
\arXiv 1806.08870

[Frob95]
F. G. Frobenius,
Verallgemeinerung des Sylow'schen Satzes,
Sitzungsberichte der K\"onigl. Preu\SS. Akad. der Wissenschaften (Berlin)
(1895), 981-993.

[Frob03]
F. G. Frobenius, \"Uber einen Fundamentalsatz der Gruppentheorie,
Sitzungsberichte der K\"onigl. Preu\SS. Akad. der Wissenschaften (Berlin)
(1903), 987-991.

[GRV12]
C. Gordon, F. Rodriguez-Villegas,
On the divisibility of $\#\Hom(\Gamma, G)$ by $|G|$,
J. Algebra, 350:1 (2012), \hbox{300-307}.
\arXiv:1105.6066

[Hall36]
Ph. Hall,
On a theorem of Frobenius,
Proc. London Math. Soc. 40 (1936), 468-501.

[Isaa70]
I. M. Isaacs,
Systems of equations and generalized characters in groups,
Canad. J. Math., 22 (1970), \hbox{1040-1046}.

[Iwa82]
S. Iwasaki,
A note on the $n$th roots ratio of a subgroup of a finite group,
J. Algebra, 78:2 (1982), 460-474.

[KM14]
A. A. Klyachko, A. A. Mkrtchyan,
How many tuples of group elements have a given property?
With an appendix by Dmitrii V. Trushin,
Intern. J. of Algebra and Comp. 24:4 (2014), 413-428.
\arXiv 1205.2824

[KM17]
A. A. Klyachko, A. A. Mkrtchyan,
Strange divisibility in groups and rings,
Arch. Math. 108:5 (2017), 441-451.
\arXiv 1506.08967

[KR20]
A. A. Klyachko, M. A. Ryabtseva,
The dimension of solution sets to systems of equations in algebraic groups,
Israel Journal of Mathematics, 237:1 (2020), 141-154.
\arXiv 1903.05236

[Kula38]
A. Kulakoff,
Einige Bemerkungen zur Arbeit: ``On a theorem of Frobenius" von P. Hall,
Mat. Sb.,  3(45):2 (1938), 403-405.

[SaAs07]
J. Sato, T. Asai,
On the $n$-th roots of a double coset of a finite group,
J. School Sci. Eng., Kinki Univ., 43 (2007), 1-4.

[Sehg62]
S. K. Sehgal,
On P. Hall's generalisation of a theorem of Frobenius,
Proc. Glasgow Math. Assoc.,  5 (1962), 97-100.

[Solo69]
L. Solomon,
The solution of equations in groups,
Arch. Math., 20:3 (1969), 241-247.


[Stru95]
S. P. Strunkov,
On the theory of equations in finite groups,
Izvestiya: Mathematics, 59:6 (1995), 1273-1282.

[Yosh93]
T. Yoshida,
$|\Hom(A, G)|$,
Journal of Algebra, 156:1 (1993), 125-156.

\end